\documentclass{amsart}%
\usepackage{amsmath}
\usepackage{color}
\usepackage{amssymb}
\usepackage{amsfonts}
\usepackage{graphicx}%
\newtheorem{theorem}{Theorem}[section]
\theoremstyle{plain}

\newtheorem{problem}{Problem}
\newtheorem{proposition}[theorem]{Proposition}

\numberwithin{equation}{section}

\begin{document}
\title[A supercritical elliptic problem]{Symmetries, Hopf fibrations and supercritical elliptic problems}
\author{M\'{o}nica Clapp}
\address{Instituto de Matem\'{a}ticas, Universidad Nacional Aut\'{o}noma de M\'{e}xico,
Circuito Exterior, C.U., 04510 M\'{e}xico D.F., Mexico}
\email{monica.clapp@im.unam.mx}
\author{Angela Pistoia}
\address{Dipartimento di Metodi e Modelli Matematici, Universit\'{a} di Roma "La
Sapienza", via Antonio Scarpa 16, 00161 Roma, Italy}
\email{pistoia@dmmm.uniroma1.it}
\thanks{Research supported by CONACYT grant 129847 and PAPIIT grant IN106612 (Mexico)
and Universit\`{a} degli Studi di Roma "La Sapienza" Accordi Bilaterali
"Esistenza e propriet\`{a} geometriche di soluzioni di equazioni ellittiche
non lineari" (Italy).}
\date{\today}
\maketitle

\begin{abstract}
We consider the semilinear elliptic boundary value problem
\[
-\Delta u=\left\vert u\right\vert ^{p-2}u\text{ in }\Omega,\text{\quad
}u=0\text{ on }\partial\Omega,
\]
in a bounded smooth domain $\Omega$ of $\mathbb{R}^{N}$ for supercritical
exponents $p>\frac{2N}{N-2}.$

Until recently, only few existence results were known. An approach which has
been successfully applied to study this problem, consists in reducing it to a
more general critical or subcritical problem, either by considering rotational
symmetries, or by means of maps which preserve the Laplace operator, or by a
combination of both.

The aim of this paper is to illustrate this approach by presenting a selection
of recent results where it is used to establish existence and multiplicity or
to study the concentration behavior of solutions at supercritical exponents.

\medskip


\end{abstract}

\section{Introduction}

Consider the model problem%
\begin{equation}
\left\{
\begin{array}
[c]{ll}%
-\Delta u=\left\vert u\right\vert ^{p-2}u & \text{in }\Omega,\\
\text{ \ \ \ \ }u=0 & \text{on }\partial\Omega,
\end{array}
\right.  \tag{$\wp_p$}\label{mainprob}%
\end{equation}
where $\Delta$ is the Laplace operator, $\Omega$ is a bounded domain in
$\mathbb{R}^{N}$ with smooth boundary, $N\geq3,$ and $p>2.$

Despite its simple form, this problem has been an amazing source of open
problems, and the process of understanding it has helped develop new and
interesting techniques which can be applied to a wide variety of problems.

The behavior of this problem depends strongly on the exponent $p$. It is
called subcritical, critical or supercritical depending on whether
$p\in(2,2_{N}^{\ast}),$ $p=2_{N}^{\ast}$ or $p\in(2_{N}^{\ast},\infty),$ where
$2_{N}^{\ast}:=\frac{2N}{N-2}$ is the so-called critical Sobolev exponent.

In the subcritical case, standard variational methods yield the existence of a
positive solution and infinitely many sign changing solutions. But if $p$ is
critical or supercritical the existence of solutions becomes a delicate issue.
It depends on the domain. An identity obtained by Pohozhaev \cite{po} implies
that (\ref{mainprob}) does not have a nontrivial solution if $\Omega$ is
strictly starshaped and $p\in\lbrack2_{N}^{\ast},\infty)$. On the other hand,
Kazdan and Warner \cite{kw} showed that infinitely many radial solutions exist
for every $p\in(2,\infty)$ if $\Omega$ is an annulus.

The critical problem has received much attention during the last thirty years,
partly due to the fact that it is a simple model for equations which arise in
some fundamental questions in differential geometry, like the Yamabe problem
or the prescribed scalar curvature problem. Still, many questions remain open
in this case.

Until quite recently, only few existence results were known for $p\in
(2_{N}^{\ast},\infty)$. A fruitful approach which has been applied in recent
years to treat supercritical problems consists in reducing problem
(\ref{mainprob}) to a more general elliptic critical or subcritical problem,
either by considering rotational symmetries, or by means of maps which
preserve the Laplace operator, or by a combination of both.

The aim of this paper is to illustrate this approach by presenting a selection
of recent results where it is used to establish existence and multiplicity or
to study the concentration behavior of solutions at certain supercritical
exponents. To put these results into perspective, we first present some
nonexistence results.

\section{Nonexistence results}

When $p=2_{N}^{\ast}$ a remarkable result obtained by Bahri and Coron
\cite{bc} establishes the existence of at least one positive solution to
problem (\ref{mainprob}) in every domain $\Omega$ having nontrivial reduced
homology with $\mathbb{Z}/2$-coefficients. Moreover, if $\Omega$ is invariant
under the action of a closed subgroup $G$ of the group $O(N)$ of linear
isometries of $\mathbb{R}^{N}$ and every $G$-orbit in $\Omega$ is
infinite\footnote{Recall that the $G$-orbit of a point $x\in\mathbb{R}^{n}$ is
the set $Gx:=\{gx:g\in G\}.$ A subset $X$ of $\mathbb{R}^{n}$ is said to be
$G$-invariant if $Gx\subset X$ for every $x\in X,$ and a function
$f:X\rightarrow\mathbb{R}$ is $G$-invariant if it is constant on every
$G$-orbit of $X.$}, the critical problem is known to have infinitely many
solutions \cite{c}.

Passaseo showed in \cite{p1,p2}\ that neither of these conditions is enough to
guarantee existence in the supercritical case. He proved the following result.

\begin{theorem}
\label{thm:passaseo}For each $1\leq k\leq N-3$ there is a domain $\Omega$ such that

\begin{enumerate}
\item[(a)] $\Omega$ has the homotopy type of $\mathbb{S}^{k},$

\item[(b)] $\Omega$ is $O(k+1)$-invariant with infinite $O(k+1)$-orbits,

\item[(c)] $(\wp_{p})$ has no solution for $p\geq2_{N,k}^{\ast}:=\frac
{2(N-k)}{N-k-2},$

\item[(d)] $(\wp_{p})$ has infinitely many solutions for $p<2_{N,k}^{\ast}.$
\end{enumerate}
\end{theorem}

Here $O(k+1)$ is the group of all linear isometries of $\mathbb{R}^{k+1}$
acting on the first $k+1$ coordinates of a point in $\mathbb{R}^{N}.$ Note
that \
\[
2_{N,k}^{\ast}:=\frac{2(N-k)}{N-k-2}=2_{N-k}^{\ast}%
\]
is the critical Sobolev exponent in dimension $N-k.$ It is called the
$(k+1)$-\emph{st critical exponent in dimension }$N.$ Note also that
$2_{N}^{\ast}<2_{N,1}^{\ast}<2_{N,2}^{\ast}<\cdots<2_{N,N-3}^{\ast}=6,$ so
$2_{N,k}^{\ast}$ is supercritical in dimension $N$ if $k\geq1.$

Passaseo's domains are defined as%
\[
\Omega:=\{(y,z)\in\mathbb{R}^{k+1}\times\mathbb{R}^{N-k-1}:\left(  \left\vert
y\right\vert ,z\right)  \in B\},
\]
where $B$ is any open ball, centered in $\left(  0,\infty\right)
\times\{0\},$ whose closure is contained in the halfspace $\left(
0,\infty\right)  \times\mathbb{R}^{N-k-1}.$%
\[%
{\includegraphics[
height=1.4624in,
width=1.7711in
]%
{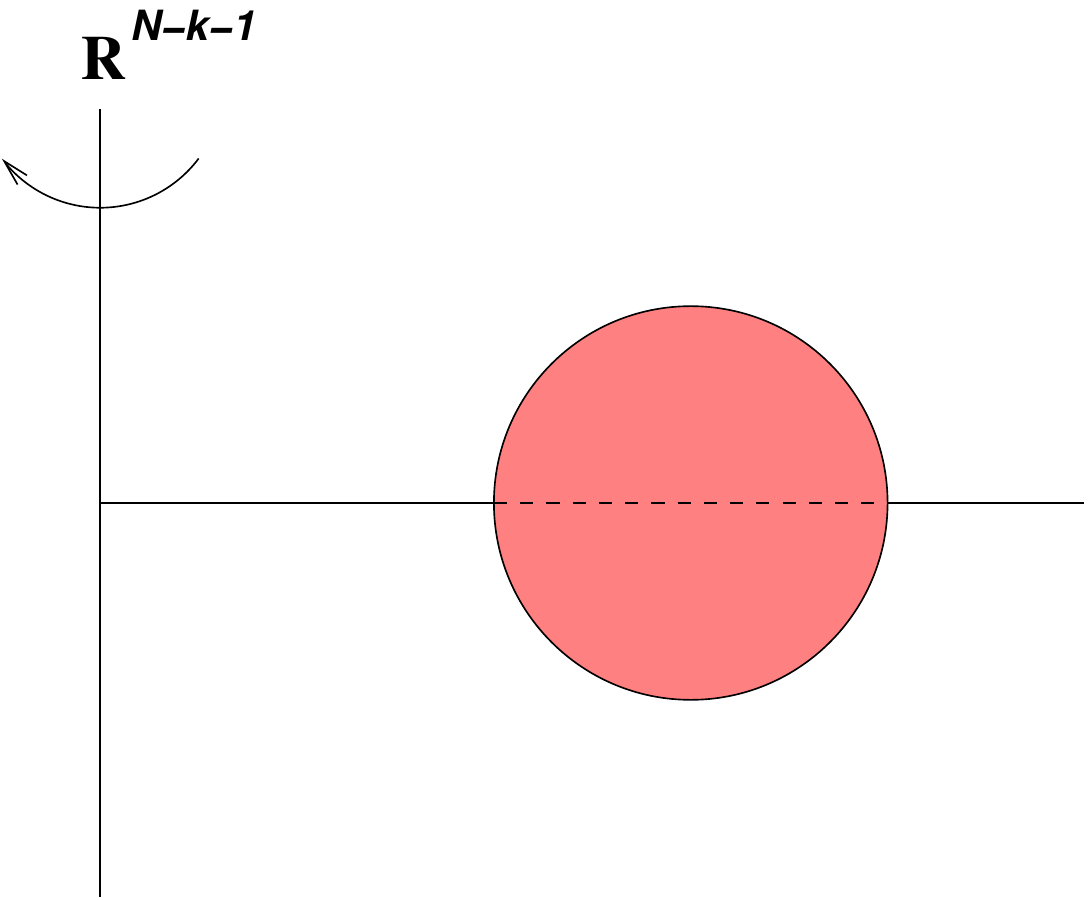}%
}
\]

Passaseo's result was extended to more general domains by Faya and the authors
of this paper in \cite{cfp1}. They also showed that existence may fail even in
domains with richer topology. More precisely, they proved the following result.

\begin{theorem}
\label{thm:cfp_nonexistence}For every $\varepsilon>0$ there is a domain
$\Omega$ such that

\begin{enumerate}
\item[(a)] $\Omega$ has the homotopy type of $\mathbb{S}^{1}\times\cdots
\times\mathbb{S}^{1}$ ($k$ factors),

\item[(b)] $\Omega$ is $O(2)\times\cdots\times O(2)$-invariant with infinite orbits,

\item[(c)] $(\wp_{p})$ does not have a nontrivial solution for $p\geq
2_{N,k}^{\ast}+\varepsilon,$

\item[(d)] $(\wp_{p})$ has infinitely many solutions for $p<2_{N,k}^{\ast}.$
\end{enumerate}
\end{theorem}

Here $\mathbb{S}^{1}$ stands for the unit circle in $\mathbb{R}^{2},$ and the
group $O(2)\times\cdots\times O(2)$ with $k$ factors acts on $\mathbb{S}%
^{1}\times\cdots\times\mathbb{S}^{1}$ in the obvious way. Note that there are
$k$\ cohomology classes in $H^{1}(\Omega;\mathbb{Z})$\ whose cup-product is
nonzero. In fact, the cup-length of $\Omega$ is $k+1.$

These domains are of the form%
\begin{equation}
\Omega:=\{(y^{1},\ldots,y^{k},z)\in\mathbb{R}^{2}\times\cdots\times
\mathbb{R}^{2}\times\mathbb{R}^{N-2k}:(|y^{1}|,\ldots,|y^{k}|,z)\in B\},
\label{toroidal}%
\end{equation}
where $B$ is an open ball centered in $\left(  0,\infty\right)  ^{k}%
\times\{0\},$ whose closure is contained in $\left(  0,\infty\right)
^{k}\times\mathbb{R}^{N-2k}$ and whose radius decreases as $\varepsilon
\rightarrow0.$%
\[%
{\includegraphics[
height=1.4633in,
width=1.6103in
]%
{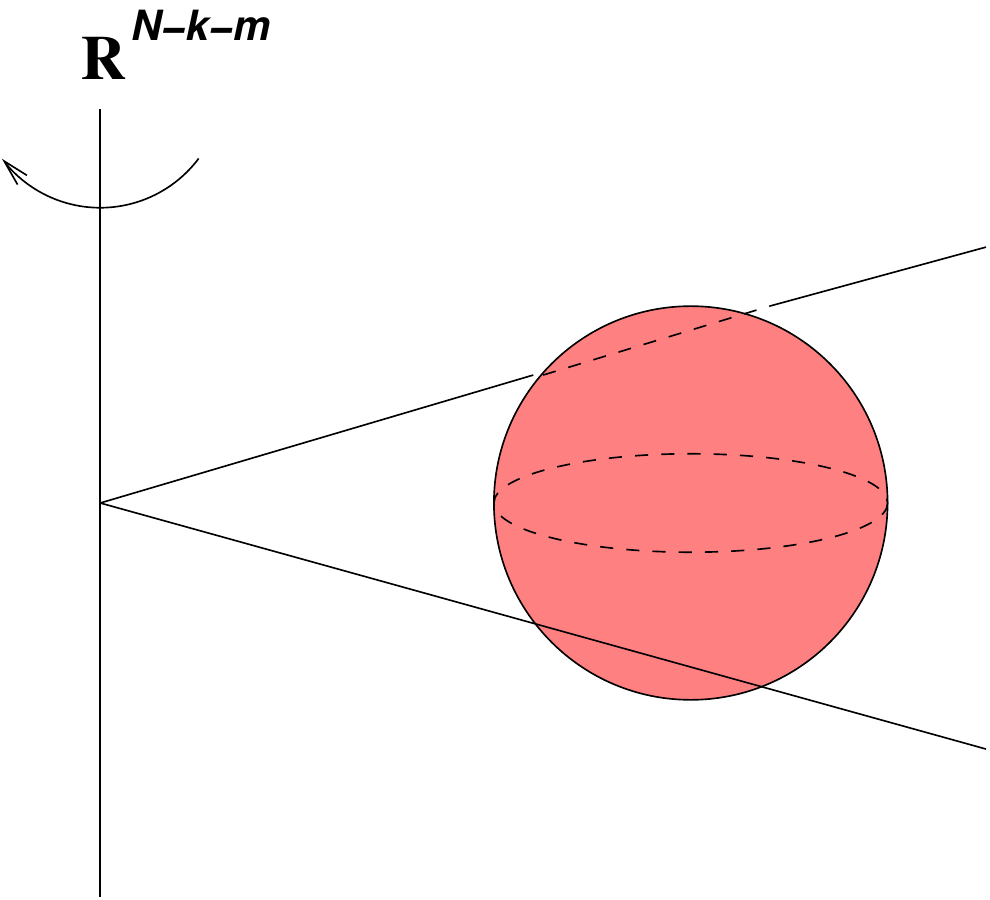}%
}
\]
The question whether one can get rid of $\varepsilon$ remains open.

\begin{problem}
Is it true that $(\wp_{p})$ does not have a nontrivial solution in the domain
$\Omega$ defined in \emph{(\ref{toroidal})} if $p\geq2_{N,k}^{\ast}$ and $B$
is an open ball of arbitrary radius?
\end{problem}

\section{Solutions for higher critical exponents via Hopf maps}

A fruitful approach to produce solutions to a supercritical problem $(\wp
_{p})$ is to reduce it to some problem of the form%
\begin{equation}
-\text{div}(a(x)\nabla v)=b(x)|v|^{p-2}v\text{ \ in}\ \Theta,\text{\quad
}v=0\text{ \ on}\ \partial\Theta, \label{general_reduction}%
\end{equation}
in a bounded smooth domain $\Theta$ in $\mathbb{R}^{n},$ with $n:=\dim
\Theta<\dim\Omega=N$ and the same exponent $p$. Thus, if $p\in(2_{N}^{\ast
},2_{n}^{\ast}],$ then $p$ is subcritical or critical for problem
(\ref{general_reduction}) but it is supercritical for (\ref{mainprob}).

\subsection{A reduction via Hopf maps}

Hopf maps provide a way to obtain such a reduction.

For $N=2,4,8,16$ we write $\mathbb{R}^{N}\mathbb{=K}\times\mathbb{K}$, where
$\mathbb{K}$ is either the real numbers $\mathbb{R}$, or the complex numbers
$\mathbb{C}$, or the quaternions $\mathbb{H},$ or the Cayley numbers
$\mathbb{O}.$

The Hopf map $\mathfrak{h}_{\mathbb{K}}:\mathbb{R}^{N}=\mathbb{K}%
\times\mathbb{K}\rightarrow\mathbb{R}\times\mathbb{K}=\mathbb{R}^{(N/2)+1}$ is
given by%
\[
\mathfrak{h}_{\mathbb{K}}(z_{1},z_{2})=(\left\vert z_{1}\right\vert
^{2}-\left\vert z_{2}\right\vert ^{2},\,2\overline{z_{1}}z_{2})\text{.}%
\]
Topologically, it is just the quotient map of $\mathbb{K}\times\mathbb{K}$
onto its orbit space under the action of $\mathbb{S}_{\mathbb{K}}:=\{\zeta
\in\mathbb{K}:\left\vert \zeta\right\vert =1\}$ given by multiplication on
each coordinate, i.e. $\zeta(y,z):=(\zeta y,\zeta z)$ for $\zeta\in
\mathbb{S}_{\mathbb{K}},$ $(y,z)\in\mathbb{K}\times\mathbb{K}$.

But, regarding our problem, the most relevant property of $\mathfrak{h}%
_{\mathbb{K}}$ is of geometric nature. It is the fact that $\mathfrak{h}%
_{\mathbb{K}}$ preserves the Laplace operator. Maps with this property are
called harmonic morphisms \cite{bw,w}. The following statement can be derived
by straightforward computation or from the general theory of harmonic morphisms.

\begin{proposition}
\label{hopf}Let $\Omega\subset\mathbb{K}^{2}$ be an $\mathbb{S}_{\mathbb{K}}%
$-invariant domain such that $0\notin\overline{\Omega}$. Set $\Theta
:=\mathfrak{h}_{\mathbb{K}}(\Omega).$ Then $u$ is an $\mathbb{S}_{\mathbb{K}}%
$-invariant solution to problem $(\wp_{p})$ iff the function $v$ given by
$u=v\circ\mathfrak{h}_{\mathbb{K}}$ is a solution to problem%
\begin{equation}
-\Delta v=\frac{1}{2\left\vert x\right\vert }|v|^{p-2}v\quad\text{in}%
\ \Theta,\qquad v=0\quad\text{on}\ \partial\Theta. \label{hopfprob}%
\end{equation}

\end{proposition}

Note that, if $\mathbb{K=C}$, $\mathbb{H}$ or $\mathbb{O}$, then $\dim
\Theta=\dim\mathbb{K}+1<2\dim\mathbb{K}=\dim\Omega.$ Therefore, $p:=2_{\dim
\mathbb{K}+1}^{\ast}=2_{N,\dim\mathbb{K}-1}^{\ast}$ is critical for
(\ref{hopfprob}) and supercritical for (\ref{mainprob}).

Recently Pacella and Srikanth showed that the real Hopf map provides a
one-to-one correspondence between $\left[  O(m)\times O(m)\right]  $-invariant
solutions of (\ref{mainprob}) in a domain $\Omega$ in $\mathbb{R}^{2m}$ and
$O(m)$-invariant solutions of (\ref{hopfprob}) in some domain $\Theta$ in
$\mathbb{R}^{m+1}$, where $O(m)$ acts on the last $m$ coordinates of
$\mathbb{R}^{m+1}\equiv\mathbb{R}\times\mathbb{R}^{m}.$ In \cite{ps} they
proved the following result.

\begin{proposition}
\label{ps}Let $N=2m$ and $\Omega$ be an $\left[  O(m)\times O(m)\right]
$-invariant bounded smooth domain in $\mathbb{R}^{2m}$ such that
$0\notin\overline{\Omega}.$ Set
\[
\Theta:=\{(t,\zeta)\in\mathbb{R}\times\mathbb{R}^{m}:\mathfrak{h}_{\mathbb{R}%
}(\left\vert y_{1}\right\vert ,\left\vert y_{2}\right\vert )=(t,\left\vert
\zeta\right\vert )\text{ for some }(y_{1},y_{2})\in\Omega\}.
\]
Then $v(t,\zeta)=w(t,\left\vert \zeta\right\vert )$ is an $O(m)$-invariant
solution of problem \emph{(\ref{hopfprob})} iff $u(y_{1},y_{2})=w(\mathfrak{h}%
_{\mathbb{R}}(\left\vert y_{1}\right\vert ,\left\vert y_{2}\right\vert ))$ is
an $\left[  O(m)\times O(m)\right]  $-invariant solution of problem
\emph{(\ref{mainprob})}.
\end{proposition}

\subsection{Multiplicity results in symmetric domains}

The previous propositions suggest to study the critical problem%
\begin{equation}
\left\{
\begin{array}
[c]{ll}%
-\Delta v=b(x)\left\vert v\right\vert ^{2_{n}^{\ast}-2}v & \text{in }\Theta,\\
\text{ \ \ \ }v=0 & \text{on }\partial\Theta,
\end{array}
\right.  \tag{$\wp_b^\ast$}\label{crithopf}%
\end{equation}
in a bounded smooth domain $\Theta$ in $\mathbb{R}^{n}$, $n\geq3,$ where
$b:\overline{\Theta}\rightarrow\mathbb{R}\ $is a positive continuous function.

This problem is variational but, due to the lack of compactness of the
associated energy functional, classical variational methods cannot be applied
to establish existence of solutions.

Under suitable symmetry assumptions compactness is restored: if $G$ is a
closed subgroup of the group $O(n)$ of linear isometries of $\mathbb{R}^{n},$
$\Theta$ and $b$ are $G$-invariant, and every $G$-orbit in $\overline{\Theta}$
has infinite cardinality, problem (\ref{crithopf}) is known to have infinitely
many $G$-invariant solutions \cite{c}. This fact,\ together with Proposition
\ref{hopf}, provides examples of domains in which problem (\ref{mainprob}) has
infinitely many solutions for some higher critical exponents. For example, one
has the following result.

\begin{theorem}
Let $\Theta$ be a solid of revolution around the $z$-axis in $\mathbb{R}^{3}$
whose closure does not intersect the $z$-axis and set $\Omega:=\mathfrak{h}%
_{\mathbb{C}}^{-1}(\Theta).$ Then the supercritical problem
\[
(\wp_{2_{4,1}^{\ast}})\text{\qquad}-\Delta u=\left\vert u\right\vert
^{4}u\quad\text{in}\ \ \Omega,\qquad u=0\quad\text{on}\ \ \partial\Omega,
\]
has infinitely many solutions which are constant on $\mathfrak{h}_{\mathbb{C}%
}^{-1}\{(r\cos\vartheta,r\sin\vartheta,t):\vartheta\in\lbrack0,2\pi]\}$ for
each $(r,0,t)\in\Theta.$
\end{theorem}

Note that $\Omega$ is homeomorphic to $\Theta\times\mathbb{S}^{1}.$ Similar
results for problems $(\wp_{2_{8,3}^{\ast}})$ and $(\wp_{2_{16,7}^{\ast}})$
can be derived using the Hopf maps $\mathfrak{h}_{\mathbb{H}}$ and
$\mathfrak{h}_{\mathbb{O}}.$

On the other hand, if $\Theta$ contains a finite orbit, problem ($\wp
_{b}^{\ast}$) might not have a nontrivial solution, as occurs when $\Theta$ is
a ball centered at the origin in $\mathbb{R}^{n}$ and $b\equiv1.$ Conditions
which guarantee that problem (\ref{crithopf}) has a prescribed number of
solutions in domains having finite orbits were recently obtained by Faya and
the authors. They proved the following result in \cite{cfp1}. Special cases of
it were previusly established in \cite{cf, cp}.

\begin{theorem}
\label{cfp_hopf}Fix a closed subgroup $\Gamma$ of $O(n)$ and a bounded smooth
$\Gamma$-invariant domain $D$ in $\mathbb{R}^{n}$ such that every $\Gamma
$-orbit in $D$ has infinite cardinality. Assume that $b$ is $\Gamma
$-invariant. Then there exists a sequence of real numbers $(\ell_{m})$ with
the following property: if $\Theta\supset D$ and $\Theta$ is invariant under
the action of a subgroup $G$ of $\Gamma$ such that%
\begin{equation}
\min_{x\in\Theta}\frac{\#Gx}{b(x)^{\frac{n-2}{2}}}>\ell_{m}, \label{ineq}%
\end{equation}
then problem \emph{(\ref{crithopf})} has at least $m$ pairs $\pm u$ of
$G$-invariant solutions in $\Theta;$ one pair does not change sign and the
rest are sign changing.
\end{theorem}

Here $\#Gx$ denotes the cardinality of the $G$-orbit of $x.$

We illustrate this result with an example. Let $D$ be a torus of revolution
around the $z$-axis in $\mathbb{R}^{3}$, $\Gamma$ be the group of all
rotations around the $z$-axis, and $b(x):=\frac{1}{2\left\vert x\right\vert }%
$. Fix $\varepsilon>0$ smaller than the distance of $D$ to the $z$-axis. Let
$G_{k}$ be the group generated by the rotation of angle $\frac{2\pi}{k}.$ If
$k\sqrt{2\varepsilon}>\ell_{m}$ and $\Theta$ is a $G_{k}$-invariant domain
which contains $D$ whose distance to the $z$-axis is at least $\varepsilon,$
then every $G_{k}$-orbit in $\Theta$ has cardinality $k$ and the inequality
(\ref{ineq}) holds true. We call a $G_{k}$-invariant domain which contains $D$
and does not intersect the\ $z$-axis a $k$-\emph{teething toy}.%
\[%
{\includegraphics[
natheight=1.458100in,
natwidth=1.458100in,
height=1.4935in,
width=1.4935in
]%
{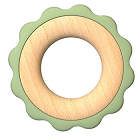}%
}
\]
Applying Proposition \ref{hopf}\ to this example we obtain the following result.

\begin{theorem}
If $\Theta$ is a $k$-teething toy whose distance to the $z$-axis is at least
$\varepsilon$ and $k\sqrt{2\varepsilon}>\ell_{m}$, then the supercritical
problem%
\[
(\wp_{2_{4,1}^{\ast}})\text{\qquad}-\Delta u=\left\vert u\right\vert
^{4}u\quad\text{in}\ \ \Omega,\qquad u=0\quad\text{on}\ \ \partial\Omega,
\]
has $m$ pairs of $\mathbb{S}_{\mathbb{C}}$-invariant solutions in
$\Omega:=\mathfrak{h}_{\mathbb{C}}^{-1}(\Theta)$; one of them does not change
sign and the rest are sign changing.
\end{theorem}

A more general result can be derived from Theorem \ref{cfp_hopf}, as stated in
\cite{cfp1}.

\subsection{Existence in domains with thin spherical holes}

Next we consider problem (\ref{crithopf}) in a punctured domain%
\[
\Theta_{\varepsilon}:=\{x\in\Theta:\left\vert x-\xi\right\vert >\varepsilon
\},
\]
where $\Theta$ is a bounded smooth domain in $\mathbb{R}^{n}$, $n\geq3,$
$\xi\in\Theta,$ $\varepsilon>0$ is small, and $b:\overline{\Theta}%
\rightarrow\mathbb{R}\ $is a positive $\mathcal{C}^{2}$-function.
Additionally, we assume that $\Theta$ and $b$ are invariant under the action
of some closed subgroup $G$ of $O(n)$ and that $\xi$ is a fixed point, i.e.
$g\xi=\xi$ for all $g\in G.$

If $b\equiv1$ Coron showed in \cite{co}\ that problem $(\wp_{2_{n}^{\ast}})$
has a positive solution in $\Theta_{\varepsilon}$ for $\varepsilon$ small
enough. Coron's proof takes advantage of the fact that the variational
functional associated to problem $(\wp_{2_{n}^{\ast}})$\ satisfies the
Palais-Smale condition between the ground state level and twice that level.
This is not true anymore when $b\not \equiv 1,$ so Coron's argument does not
carry over to the nonautonomous case.

A method which has proved to be very successful in dealing with critical
problems which involve small perturbations of the domain is the
Lyapunov-Schmidt reduction method, see e.g. \cite{gmp} and the references
therein. This method was used by Faya and the authors in \cite{cfp2} to prove
the following result.

\begin{theorem}
\label{thm:cfp_punctured}If $\nabla b(\xi)\neq0$ then, for $\varepsilon$ small
enough,
\[
(\wp_{b}^{\ast})\text{\qquad}-\Delta v=b(x)\left\vert v\right\vert
^{2_{n}^{\ast}-2}v\quad\text{in \ }\Theta_{\varepsilon},\qquad v=0\quad
\text{on \ }\partial\Theta_{\varepsilon},
\]
has a positive $G$-invariant solution $v_{\varepsilon}$ in $\Theta
_{\varepsilon}$ which concentrates at the boundary of the hole and blows up at
$\xi$ as $\varepsilon\rightarrow0$.
\end{theorem}

Note that $G$ may be the trivial group, so this result is true in a
non-symmetric setting and, combined with Proposition \ref{hopf}, yields
solutions to supercritical problems concentrating around a spherical hole, see
\cite{cfp2}.

But we may also combine it with Proposition \ref{ps}\ as follows: Let $N=2m,$
$m\geq2,\ \Omega$ be an $\left[  O(m)\times O(m)\right]  $-invariant bounded
smooth domain in $\mathbb{R}^{N}\equiv\mathbb{R}^{m}\times\mathbb{R}^{m}$ such
that $0\notin\overline{\Omega},$ and $\xi\in\Omega\cap\left(  \mathbb{R}%
^{m}\times\{0\}\right)  .$ For $\varepsilon>0$ small enough set%
\[
\Omega_{\varepsilon}:=\{x\in\Omega:\text{dist}(x,S_{\xi})>\varepsilon\},
\]
where $S_{\xi}:=\{(x,0)\in\mathbb{R}^{m}\times\{0\}:\left\vert x\right\vert
=\left\vert \xi\right\vert \}.$ The following result, obtained by combining
Theorem \ref{thm:cfp_punctured} with Proposition \ref{ps}, was established in
\cite{cfp2}.

\begin{theorem}
For each $\varepsilon$ small enough the supercritical problem
\[
(\wp_{2_{2m,m-1}^{\ast}})\qquad-\Delta u=\left\vert u\right\vert
^{4/(m-1)}u\quad\text{in \ }\Omega_{\varepsilon},\qquad v=0\quad\text{on
\ }\partial\Omega_{\varepsilon},
\]
has a positive $\left[  O(m)\times O(m)\right]  $-invariant solution
$u_{\varepsilon}$ which concentrates along the set $\{x\in\Omega:$
\emph{dist}$(x,S_{\xi})=\varepsilon\}$ and blows up at the $(m-1)$-dimensional
sphere $S_{\xi}$ as $\varepsilon\rightarrow0$.
\end{theorem}

\section{Solutions for higher critical exponents via rotations}

Supercritical problems in domains obtained through rotations can be reduced to
subcritical or critical problems as follows.

\subsection{A reduction via rotations}

\label{subsec:rot}Fix $k_{1},\ldots,k_{m}\in\mathbb{N}$ and set $k:=k_{1}%
+\cdots+k_{m}.$ If $N\geq k+m$ let
\[
\Omega:=\{(y^{1},\ldots,y^{m},z)\in\mathbb{R}^{k_{1}+1}\times\cdots
\times\mathbb{R}^{k_{m}+1}\times\mathbb{R}^{N-k-m}:\left(  \left\vert
y^{1}\right\vert ,\ldots,\left\vert y^{m}\right\vert ,z\right)  \in\Theta\},
\]
where $\Theta$ is a bounded smooth domain in $\mathbb{R}^{N-k}$ whose closure
is contained in $(0,\infty)^{m}\times\mathbb{R}^{N-k-m}.$
\[%
{\includegraphics[
height=1.4875in,
width=2.2563in
]%
{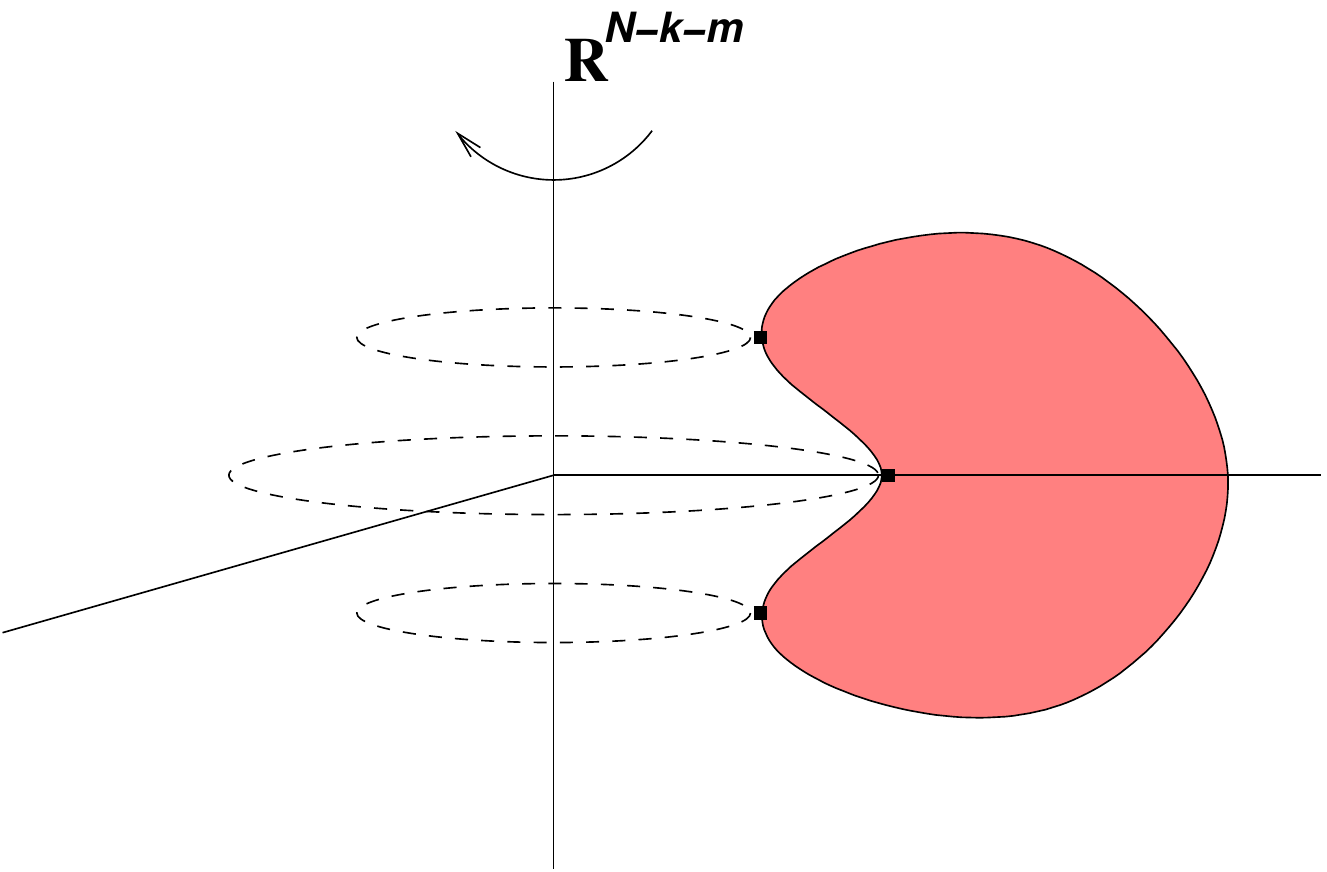}%
}
\]
Each point $\xi\in\Theta$ gives rise to a subset
\begin{equation}
T_{\xi}:=\{(y^{1},\ldots,y^{m},z)\in\Omega:\left\vert y^{i}\right\vert
=\xi_{i},\text{ }z=(\xi_{m+1},\ldots,\xi_{N-k})\} \label{txi}%
\end{equation}
of $\Omega$ which is homeomorphic to the product of spheres $\mathbb{S}%
^{k_{1}}\times\cdots\times\mathbb{S}^{k_{m}}.$ A straightforward computation
yields the following result.

\begin{proposition}
\label{rot}A function $u$ of the form $u(y^{1},\ldots,y^{m},z)=v(\left\vert
y^{1}\right\vert ,\ldots,\left\vert y^{m}\right\vert ,z)\ $is a solution of
problem $(\wp_{p})$ iff $v\ $is a solution of%
\begin{equation}
-\text{\emph{div}}(a(x)\nabla v)=a(x)|v|^{p-2}v\text{ in}\ \Theta,\text{\quad
}v=0\text{ on}\ \partial\Theta, \label{rotprob}%
\end{equation}
with $a(x_{1},\ldots,x_{N-k}):=x_{1}^{k_{1}}\cdots x_{m}^{k_{m}}.$
\end{proposition}

\subsection{Multiplicity results in domains obtained by rotation}

Wei and Yan considered domains as above, with $m=1$, where $\Theta$ is
invariant under the action of the group $O(2)\times O(1)^{N-k-3}$ on the last
$N-k-1$ coordinates of $\mathbb{R}^{N-k},$ i.e.%
\[%
\begin{array}
[c]{l}%
(s,r\cos\theta,r\sin\theta,x_{3},\ldots,x_{N-k})\in\Theta\text{ for all
}\theta\in(0,2\pi)\text{ \ if }(s,r,0,\ldots,x_{N-k})\in\Theta,\\
(x_{1},,\ldots,-x_{i},,\ldots,x_{N-k})\in\Theta\text{ \ if }(x_{1}%
,\ldots,x_{i},\ldots,x_{N-k})\in\Theta\text{ and }i=4,\ldots,N-k.
\end{array}
\]
They proved the following result in \cite{wy}.

\begin{theorem}
Let $N\geq5.$ Assume that $\Theta$ is $\left[  O(2)\times O(1)^{N-k-3}\right]
$-invariant and that there is a point
\[
(s^{\ast},r^{\ast})\in\mathcal{S}:=\{(s,r)\in\mathbb{R}^{2}:(s,r,0,\ldots
,0)\in\partial\Theta\},
\]
which is a strict local minimum (or a strict local maximum) of the distance of
$\mathcal{S}$ to $\{0\}\times\mathbb{R}.$ Set $\xi_{j}:=(s^{\ast},r^{\ast}%
\cos\left(  2\pi j/\ell\right)  ,r^{\ast}\sin\left(  2\pi j/\ell\right)
,0,\ldots,0).$ Then, for large enough $\ell\in\mathbb{N}$, problem
$(\wp_{2_{N,k}^{\ast}})$ has a solution $u_{\ell}$ with $\ell$ positive
layers; which concentrate along each of the $k$-dimensional spheres
$T_{\xi_{j}}\subset\partial\Omega,$ $j=0,\ldots,\ell-1,$ and blow up at
$\partial\Omega$ as $\ell\rightarrow\infty.$
\end{theorem}

They derived this result from Proposition \ref{rot}\ after proving the
existence of $\ell$-multibubble solutions for the critical problem%
\[
(\wp_{a,a}^{\ast})\qquad-\text{div}(a(x)\nabla v)=a(x)|v|^{2_{N-k}^{\ast}%
-2}v\text{\quad in}\ \ \Theta,\text{\qquad}v=0\text{\quad on}\ \ \partial
\Theta,
\]
which concentrate at the points $\xi_{j}\in\partial\Theta.$ A precise
description of the solutions is given in \cite{wy}.

In \cite{kp1} Kim and Pistoia considered domains with thin $k$-dimensional
holes. More precisely, for some fixed $\xi\in\Theta,$ they considered
\[
\Omega_{\varepsilon}:=\{x\in\Omega:\text{dist}(x,T_{\xi})>\varepsilon\},
\]
with $T_{\xi}$ as in (\ref{txi}) and $\varepsilon>0$ sufficiently small. Note
that $\Omega_{\varepsilon}$ is obtained by rotating the punctured domain
$\Theta_{\varepsilon}:=\{x\in\Theta:\left\vert x-\xi\right\vert >\varepsilon
\}$ as described in subsection \ref{subsec:rot}. For $N-k\geq4$ they proved
the existence of towers of bubbles with alternating signs around $\xi$ for the
anisotropic problem $(\wp_{a,a}^{\ast})$ in $\Theta_{\varepsilon},$ thus
extending a previous result by Ge, Musso and Pistoia \cite{gmp} for the
autonomous problem $a\equiv1$. The number of bubbles increases as
$\varepsilon\rightarrow0.$ Combining this result with Proposition \ref{rot}
they obtained the following one.

\begin{theorem}
Let $N\geq k+4.$ Then, for every $\ell\in\mathbb{N}$ there exists
$\varepsilon_{\ell}>0$ such that, for each $\varepsilon\in(0,\varepsilon
_{\ell}),$ the supercritical problem%
\[
-\Delta u=\left\vert u\right\vert ^{2_{N,k}^{\ast}-2}u\quad\text{in \ }%
\Omega_{\varepsilon},\qquad v=0\quad\text{on \ }\partial\Omega_{\varepsilon},
\]
has a solution $u_{\varepsilon}$ with $\ell$ layers of alternating signs,
which concentrate with different rates along the boundary of the tubular
neighborhood of radius $\varepsilon$ of $T_{\xi}$ and blow up at $T_{\xi}$ as
$\varepsilon\rightarrow0.$
\end{theorem}

A precise description of the solutions can be found in \cite{kp1}.

\section{Concentration along manifolds at the higher critical exponents}

\subsection{Approaching the higher critical exponents from below}

For domains as those described in subsection \ref{subsec:rot} and
$p\in(2,2_{N,k}^{\ast})$ problem $(\wp_{p})$\ has infinitely many nontrivial
solutions of the form $u(y^{1},\ldots,y^{m},z)=v(\left\vert y^{1}\right\vert
,\ldots,\left\vert y^{m}\right\vert ,z).$ This follows immediately from
Proposition \ref{rot} using standard variational methods because problem
(\ref{rotprob}) is subcritical and the domain $\Theta$ is bounded. Existence
and nonexistence results in some unbounded domains are also available
\cite{cs}.

So the question is whether one can establish existence of solutions to
$(\wp_{p})$ which exhibit a certain concentration behavior as $p\rightarrow
2_{N,k}^{\ast}.$

For the slightly subcritical problem $(\wp_{2_{N}^{\ast}-\varepsilon})$
positive and sign changing solutions $u_{\varepsilon}$ which blow up at one or
several points in $\Omega$ as $\varepsilon\rightarrow0$ have been obtained
e.g. in \cite{bmp, pw, r}. Recently, del Pino, Musso and Pacard \cite{dmp}
considered the case in which $p$ approaches the second critical exponent from
below. They showed that, if $N\geq8$ and $\partial\Omega$ contains a
nondegenerate closed geodesic $\Gamma$ with negative inner normal curvature
then, for every $\varepsilon>0$ small enough, away from an explicit discrete
set of values, problem $(\wp_{2_{N,1}^{\ast}-\varepsilon})$ has a positive
solution $u_{\varepsilon}$ which concentrates and blows up at $\Gamma$ as
$\varepsilon\rightarrow0.$

It is natural to ask whether similar concentration phenomena can be observed
as $p$ aproaches the $(k+1)$-st critical exponent $2_{N,k}^{\ast}$ from below,
i.e. whether there are domains in which problem $(\wp_{2_{N,k}^{\ast
}-\varepsilon})$ has a solution $u_{\varepsilon}$ which concentrates and blows
up at a $k$-dimensional submanifold of $\Omega$ as $\varepsilon\rightarrow0.$
Ackermann, Kim and the authors have given positive answers to this question
for domains $\Omega$ as in subsection \ref{subsec:rot}.

Let $\mathcal{K}$ be the set of all nondegenerate critical points $\xi$ of the
restriction of the function $a(x_{1},\ldots,x_{N-k}):=x_{1}^{k_{1}}\cdots
x_{m}^{k_{m}}$ to $\partial\Theta$ such that $\nabla a(\xi)$ points into the
interior of $\Theta,$ and let $T_{\xi}$ be the set defined in (\ref{txi}). The
following result was proved in \cite{acp}.

\begin{theorem}
\label{thm:conc1}For any subset $\{\xi_{1},\dots,\xi_{\ell}\}$ of
$\mathcal{K}$ and $1\leq m\leq\ell$ there exists $\varepsilon_{0}>0$ such
that, for each $\varepsilon\in(0,\varepsilon_{0}),$ problem $(\wp
_{2_{N,k}^{\ast}-\varepsilon})$ has a solution $u_{\varepsilon}$ with $m$
positive layers and $\ell-m$ negative layers; which concentrate at the same
rate and blows up along one of the sets $T_{\xi_{i}}$ as $\varepsilon
\rightarrow0.$
\end{theorem}

Sign changing solutions are also available. Statement (a) in the following
theorem was proved in \cite{acp} and statement (b) was proved in \cite{kp2}.

\begin{theorem}
\label{thm:conc2}Assume there exist $\xi_{0}\in\mathcal{K}$ and $\tau
_{1},\dots,\tau_{N-k-1}\in\mathbb{R}^{N-k}$ such that the set $\{\nabla
a(\xi_{0}),\tau_{1},\dots,\tau_{N-k-1}\}$ is orthogonal and $\Theta$ and $a$
are invariant with respect to the reflection $\varrho_{i}$ on the hyperplane
through $\xi_{0}$ which is orthogonal to $\tau_{i},$ for each $i=1,\dots
,N-k-1$. Then the following statements hold true:

\begin{enumerate}
\item[(a)] For each $\varepsilon>0$ small enough problem $(\wp_{2_{N,k}^{\ast
}-\varepsilon})$ has a sign changing solution $u_{\varepsilon}$ with one
positive and one negative layer, which concentrate at the same rate along
$T_{\xi_{0}}$ and blow up at $T_{\xi_{0}}$ as $\varepsilon\rightarrow0.$

\item[(b)] If $k\leq N-4$ then, for any integer $\ell\geq2,$ there exists
$\varepsilon_{\ell}>0$ such that, for $\varepsilon\in(0,\varepsilon_{\ell}),$
problem $(\wp_{2_{N,k}^{\ast}-\varepsilon})$ has a solution $u_{\varepsilon}$
with $\ell$ layers of alternating signs which concentrate at different rates
along $T_{\xi_{0}}$ and blow up at $T_{\xi_{0}}$ as $\varepsilon\rightarrow0.$
\end{enumerate}
\end{theorem}

These results follow from Proposition \ref{rot}\ once the corresponding
statements for the slightly subcritical problems
\[
-\text{div}(a(x)\nabla v)=a(x)|v|^{2_{N-k}^{\ast}-\varepsilon-2}v\text{\quad
in}\ \ \Theta,\text{\qquad}v=0\text{\quad on}\ \ \partial\Theta,
\]
have been established. This is done in \cite{acp,kp2}, where a precise
description of the solutions is given. The following questions were raised in
\cite{acp}:

\begin{problem}
Can statement \emph{(a)} in \emph{Theorem \ref{thm:conc2}}\ be improved to
establish existence of solutions with $\ell$ layers of alternating signs which
concentrate at the same rate along $T_{\xi_{0}}$ and blow up at $T_{\xi_{0}}$
as $\varepsilon\rightarrow0,$ for any $\ell\geq2$?
\end{problem}

\begin{problem}
Does \emph{Theorem \ref{thm:conc2}}\ hold true without the symmetry assumption?
\end{problem}

Recently, Pacella and Pistoia considered the case in which $N=2m$ and $\Omega$
is the annulus $\Omega:=\{x\in\mathbb{R}^{N}:0<a<\left\vert x\right\vert
<b\}$. They used Proposition \ref{ps} to establish the existence of positive
and sign changing solutions to problem $(\wp_{2_{N,k}^{\ast}-\varepsilon})$
concentrating along the $(m-1)$-dimensional spheres%
\[
S_{1}:=\{(x,0)\in\mathbb{R}^{m}\times\{0\}:\left\vert x\right\vert
=a\}\text{,\qquad}S_{2}:=\{(0,y)\in\{0\}\times\mathbb{R}^{m}:\left\vert
y\right\vert =a\}
\]
as $\varepsilon\rightarrow0.$ Their results can be found in \cite{pp}.

\subsection{Approaching the higher critical exponents from above}

For the slightly supercritical problem $(\wp_{2_{N}^{\ast}+\varepsilon}),$
$\varepsilon>0,$ in a domain $\Omega$ with nontrivial topology del Pino,
Felmer and Musso established the existence of a positive solution
$u_{\varepsilon}$ with two bubbles which concentrate at two different points
$\xi_{1},\xi_{2}\in\Omega$ as $\varepsilon\rightarrow0$ \cite{dfm}. Solutions
with more that two bubbles are also available, see e.g. \cite{bmp,dfm2,pr}.
For $\varepsilon$ sufficiently small solutions with only one bubble do not
exist \cite{berg}.

The following problem is fully open to investigation.

\begin{problem}
Are there domains $\Omega$ in which $(\wp_{2_{N,k}^{\ast}+\varepsilon})$ has
positive or sign changing solutions $u_{\varepsilon}$ which concentrate and
blow up at $k$-dimensional manifolds as $\varepsilon\rightarrow0$?
\end{problem}

\end{document}